\documentclass[12pt, leqno]{article}

\def\diam{\mathop{\rm diam}}

\newtheorem{theorem}{Theorem}
\newtheorem{lemma}[theorem]{Lemma}
\newtheorem{proposition}[theorem]{Proposition}
\newtheorem{definition}[theorem]{Definition}
\newtheorem{remark}[theorem]{Remark}

\newcommand{\begintheorem}{\addtocounter{equation}{1}\begin{theorem}}
\newcommand{\beginlemma}{\addtocounter{equation}{1}\begin{lemma}}
\newcommand{\beginproposition}{\addtocounter{equation}{1}\begin{proposition}}
\newcommand{\begindefinition}{\addtocounter{equation}{1}\begin{definition}}
\newcommand{\beginremark}{\addtocounter{equation}{1}\begin{remark}}

\begin{document}

\title{A beginner's guide to analysis on metric spaces}

\author{Stephen William Semmes \\
        Rice University \\
        Houston, Texas}

\date{}

\maketitle

\begin{abstract}
These notes deal with various kinds of distance functions and related
properties and measurements for sets and functions.
\end{abstract}

\tableofcontents

\section{Basic notions}
\label{basic notions}
\setcounter{equation}{0}

	Let $M$ be a nonempty set.  By a \emph{distance function} on
$M$ we mean a real-valued function $d(x, y)$ defined for $x, y \in M$
such that (i) $d(x, y) \ge 0$ for all $x, y \in M$, (ii) $d(x, y) = 0$
if and only if $x = y$, and (iii)
\begin{equation}
	d(y, x) = d(x, y)
\end{equation}
for all $x, y \in M$.  If moreover the triangle inequality holds,
which is to say that
\begin{equation}
	d(x, z) \le d(x, y) + d(y, z)
\end{equation}
for all $x, y, z \in M$, then the distance function $d(\cdot, \cdot)$
is said to be a \emph{metric} on $M$.  More generally, if there is a
real number $C \ge 1$ such that
\begin{equation}
	d(x, z) \le C (d(x, y) + d(y, z))
\end{equation}
for all $x, y, z \in M$, then the distance function $d(\cdot, \cdot)$
is said to be a \emph{quasimetric} on $M$, with constant $C$.

	A basic example occurs with the real numbers ${\bf R}$
equipped with the standard distance function.  Recall that if $x$ is a
real number, then the absolute value of $x$ is denoted $|x|$ and is
defined to be equal to $x$ is $x \ge 0$, and to be equal to $-x$ if $x
\le 0$.  Thus the absolute value of a real number is always a
nonnegative real number, the absolute value of a real number is equal
to $0$ if and only if the real number itself is equal to $0$, and the
absolute value of a product of two real numbers is equal to the
product of the absolute values of the two real numbers individually.
Furthermore,
\begin{equation}
	|x + y| \le |x| + |y|
\end{equation}
for all real numbers $x$, $y$, which is also referred to as a triangle
inequality.  It follows easily that the standard distance function $|x
- y|$ on the real numbers is in fact a metric.

	A distance function $d(x, y)$ on a nonempty set $M$ is said to
be an \emph{ultrametric} if
\begin{equation}
	d(x, z) \le \max (d(x, y), d(y, z))
\end{equation}
for all $x, y, z \in M$.  Thus an ultrametric is automatically a
metric.  One might consider the weaker condition that there be a real
number $C' \ge 1$ such that
\begin{equation}
\label{weakened ultrametric condition}
	d(x, z) \le C' \max (d(x, y), d(y, z))
\end{equation}
for all $x, y, z \in M$.  Of course this would imply that the distance
function is a quasimetric, and conversely this condition holds for
quasimetrics with $C'$ equal to $2$ times the quasimetric constant.
In other words, this weaker condition is equivalent to saying that the
distance function is a quasimetric, except perhaps with a slightly
different constant.

	Let $q$ be a positive real number.  If $q \ge 1$, then
\begin{equation}
	(a + b)^q \le 2^{q-1} \, (a^q + b^q)
\end{equation}
for all nonnegative real numbers $a$, $b$.  This amounts to the
convexity of the function $t^q$ on the nonnegative real numbers.  It
follows that the distance function $|x - y|^q$ on the real line is a
quasimetric with constant $2^{q-1}$ when $q \ge 1$.  More generally,
if $d(x, y)$ is a quasimetric on a nonempty set $M$, then $d(x, y)^q$
is also a quasimetric on $M$ when $q \ge 1$, where the quasimetric
constant $C$ for $d(x, y)$ is replaced with $2^{q-1} \, C^q$.

	Now suppose that $0 < q \le 1$.  In this case we have that
\begin{equation}
	(a + b)^q \le a^q + b^q
\end{equation}
for all nonnegative real numbers $a$, $b$.  Indeed, 
\begin{eqnarray}
	a + b & \le & (a^q + b^q) \max(a^{1-q}, b^{1-q})	\\
		& = & (a^q + b^q) \max(a^q, b^q)^{(1-q)/q}	
						\nonumber	\\
		& \le & (a^q + b^q) \, (a^q + b^q)^{(1-q)/q}	
						\nonumber	\\
		& = & (a^q + b^q)^{1/q}.	\nonumber
\end{eqnarray}
Thus we get in this case that the distance function $|x - y|^q$ on the
real line is a metric when $0 < q \le 1$.  If $d(x, y)$ is a
quasimetric on a nonempty set $M$ with constant $C$, then $d(x, y)^q$
is a quasimetric on $M$ with constant $C^q$ when $0 < q \le 1$, and in
particular $d(x, y)^q$ is a metric on $M$ if $d(x, y)$ is.

	For ultrametrics the situation is much more straightforward.
Namely, if $d(x, y)$ is an ultrametric on a nonempty set $M$ and $q$
is a positive real number, then $d(x, y)^q$ is also an ultrametric on
$M$.  If more generally $d(x, y)$ is a distance function on $M$ which
satisfies (\ref{weakened ultrametric condition}) for some $C' \ge 1$,
then $d(x, y)^q$ satisfies (\ref{weakened ultrametric condition}) with
$C'$ replaced with $(C')^q$.

	Let $d(x, y)$ be a distance function on a nonempty set $M$.
We say that $d(x, y)$ is \emph{almost a metric} with constant
$C \ge 1$ if
\begin{equation}
	d(w_0, w_l) \le C \sum_{j=1}^{l} d(w_j, w_{j-1})
\end{equation}
for any finite sequence of points $w_0, w_1, \ldots, w_l$ in $M$.  If
we restrict our attention to the case where $l = 1$, then this is the
same as the quasimetric condition.  If $d(x, y)$ is a metric, then it
is also almost a metric with $C = 1$, by repeating the triangle
inequality.  However, if $q$ is a positive real number such that $q >
1$, then the distance function $|x - y|^q$ on the real line is a
quasimetric but it is not almost a metric, as one can check using
simple examples.

	Suppose that $d(x, y)$ is a distance function on a nonempty
set $M$ which is almost a metric with constant $C$.  Define $\rho(x,
y)$ for $x, y \in M$ to be the infimum of the sums $\sum_{j=1}^l
d(w_j, w_{j-1})$ over all finite sequences $w_0, \ldots, w_l$ of
points in $M$ such that $w_0 = x$ and $w_l = y$.  We have that
\begin{equation}
	C^{-1} \, d(x, y) \le \rho(x, y) \le d(x, y)
\end{equation}
for all $x, y \in M$, where the first inequality follows from the
assumption that $d(x, y)$ is almost a metric with constant $C$, and
the second inequality is automatic from the definition of $\rho(x,
y)$.  It is easy to see that $\rho(x, y)$ is a distance function on
$M$, and in fact it is a metric, because the triangle inequality for
$\rho(x, y)$ is also an automatic consequence of the definition of
$\rho(x, y)$.

	If $d(x, y)$ is a distance function on a nonempty set $M$,
then we say that $d(x, y)$ is \emph{almost an ultrametric} with
constant $C' \ge 1$ if
\begin{equation}
	d(w_0, w_l) \le C' \max(d(w_0, w_1), \ldots, d(w_{l-1}, w_l))
\end{equation}
for any finite sequence of points $w_0, \ldots, w_l$ in $M$.  When $l
= 1$ this is the same as (\ref{weakened ultrametric condition}), and
an ultrametric is almost an ultrametric with constant equal to $1$.
On the real line, the distance function $|x - y|^q$ is not almost an
ultrametric for any positive real number $q$, as one can verify.

	Let $d(x, y)$ be a distance function on a nonempty set $M$
which is almost an ultrametric with constant $C'$.  Define $\sigma(x, y)$
for $x, y \in M$ to be the infimum of
\begin{equation}
	\max(d(w_0, w_1), \ldots, d(w_{l-1}, w_l))
\end{equation}
over all finite sequences $w_0, \ldots, w_l$ of points in $M$
with $w_0 = x$ and $w_l = y$.  Thus
\begin{equation}
	(C')^{-1} \, d(x, y) \le \sigma(x, y) \le d(x, y)
\end{equation}
for all $x, y \in M$, where the first inequality follows from the
assumption that $d(x, y)$ is almost an ultrametric with constant $C'$,
and the second condition is immediate from the definition.  Clearly
$\sigma(x, y)$ is a distance function, and in fact it is an
ultrametric, because it satisfies the ultrametric version of the
triangle inequality by construction.

	Suppose that $M$ is a nonempty set, $E$ is a nonempty subset
of $M$, and that $d(x, y)$ is a distance function on $M$.  One can
restrict $d(x, y)$ to $E$ to get a distance function on $E$.
If $d(x, y)$ is a metric, quasimetric, ultrametric, etc., on $M$,
then the same condition will hold automatically on $E$ as well.
However, depending on the structure of $E$, it may be that the
restriction of the distance function to $E$ satisfies conditions
that do not work on $M$.

	As a basic example, let $M$ be the real line, and let $E$ be
the usual Cantor middle-thirds set, as in \cite{Rudin}.  By
definition, one starts with the unit interval $[0, 1]$ in the real
line, consisting of all real numbers $x$ such that $0 \le x \le 1$,
one removes the open middle third of this interval, leaving the two
intervals $[0, 1/3]$ and $[2/3, 1]$, one removes the open middle
thirds of these two intervals to get four closed intervals of length
$1/9$, and so on.  At the $k$th stage of the construction one has a
set consisting of $2^k$ closed intervals of length $1/3^k$, and the
union of these intervals can be denoted $E_k$.  By construction,
$E_{k+1} \subseteq E_k$ for all $k$, and the Cantor set $E$ is defined
to be the intersection of all of the $E_k$'s.  The restriction of the
usual distance function $|x - y|$ to the Cantor set is almost an
ultrametric.

	For both the notions of a distance function $d(x, y)$ being
almost a metric or almost an ultrametric, one can consider more
general conditions with a lower bound in terms of some function of
$d(x, y)$, in place of simply a positive multiple of $d(x, y)$.  One
could still define $\rho(x, y)$ or $\sigma(x, y)$ as before, with
$\rho(x, y), \sigma(x, y) \le d(x, y)$ automatically, and with lower
bounds in terms of a function of $d(x, y)$.  One would again have that
$\rho(x, y)$ is a metric or that $\sigma(x, y)$ is an ultrametric, as
appropriate, by construction.  This type of condition occurs for
instance with more general Cantor sets, in which the sizes of the open
intervals being removed are more variable.  Of course these sizes are
also relevant for the constants involved even if one has lower bounds
in terms of multiples of the original distance function.

\beginremark
{\rm A result of Mac\'{\i}as and Segovia \cite{M-S} states that
for each positive real number $C_1 \ge 1$ there are positive
real numbers $a$, $C_2$ with $a \le 1$ and $C_2 \ge 1$
such that the following is true.  Let $M$ be a nonempty
set and $d(x, y)$ a distance function on $M$.  If $d(x, y)$
is a quasimetric with constant $C_1$, then $d(x, y)^a$
is almost a metric with constant $C_2$.
}
\end{remark}

	Let $n$ be a positive integer, and let ${\bf R}^n$ denote the
usual space of $n$-tuples of real numbers.  Thus $x \in {\bf R}^n$ can
be written explicitly as $x = (x_1, \ldots, x_n)$, where each
component $x_j$, $1 \le j \le n$, is a real number.  If $x, y \in {\bf
R}^n$, then $x + y$ is defined to be the element of ${\bf R}^n$ whose
$j$th component is equal to the sum of the $j$th components of $x$,
$y$.  If $x \in {\bf R}^n$ and $t$ is a real number, then the scalar
product $t \, x$ is defined to be the element of ${\bf R}^n$ whose
$j$th component is equal to the product of $t$ and the $j$th component
of $x$.  In this manner ${\bf R}^n$ is a real vector space of
dimension $n$.

	Let $N(x)$ be a real-valued function on ${\bf R}^n$
such that $N(x) \ge 0$ for all $x \in {\bf R}^n$, $N(x) = 0$
if and only if $x = 0$, and 
\begin{equation}
	N(t \, x) = |t| \, N(x)
\end{equation}
for all real numbers $t$ and all $x \in {\bf R}^n$.  If also
$N(x)$ satisfies the triangle inequality
\begin{equation}
	N(x + y) \le N(x) + N(y)
\end{equation}
for all $x, y \in {\bf R}^n$, then $N$ is said to define a \emph{norm}
on ${\bf R}^n$.  In the presence of the other conditions, one can
check that the triangle inequality is equivalent to the convexity of
the closed unit ball associated to $N$.  In other words, in the
presence of the other conditions, $N$ is a norm if and only if for
each $x, y \in {\bf R}^n$ with $N(x), N(y) \le 1$ and each real number
$t$ with $0 \le t \le 1$ we have that $N(t \, x + (1 - t) \, y) \le 1$
too.  As a generalization of the triangle inequality, if $N$ satisfies
the first set of conditions, then we say that $N$ is a
\emph{quasinorm} with constant $C \ge 1$ if
\begin{equation}
	N(x + y) \le C (N(x) + N(y))
\end{equation}
for all $x, y \in {\bf R}^n$.

	As a basic family of examples, let $p$ be a positive real
number, and define $\|x\|_p$ for $x \in {\bf R}^n$ by
\begin{equation}
	\|x\|_p = \biggl(\sum_{j=1}^n |x_j|^p \biggr)^{1/p}.
\end{equation}
Let us extend this to $p = \infty$ by putting
\begin{equation}
	\|x\|_\infty = \max(|x_1|, \ldots, |x_n|).
\end{equation}
These functions clearly satisfy the first set of conditions mentioned
in the previous paragraph, namely $\|x\|_p$ is a nonnegative real
number which is equal to $0$ exactly when $x = 0$, and $\|t \, x\|_p$
is equal to the product of $|t|$ and $\|x\|_p$ for all real numbers
$t$ and $x \in {\bf R}^n$.

	When $p = 1, \infty$ it is easy to see directly from the
definition that $\|x\|_p$ satisfies the triangle inequality and is
therefore a norm on ${\bf R}^n$.  This is also true when $1 < p <
\infty$, as one can see by checking the convexity of the closed unit
ball of $\|x\|_p$.  To be more precise, one can reduce this to the
convexity of the function $u^p$ on the nonnegative real numbers
when $p \ge 1$.  When $0 < p < 1$ and $n \ge 2$, $\|x\|_p$ does
not define a norm on ${\bf R}^n$.  It does define a quasinorm,
and indeed one has that
\begin{equation}
	\|x + y\|_p^p \le \|x\|_p^p + \|y\|_p^p
\end{equation}
for all $x, y \in {\bf R}^n$ when $0 < p < 1$.

	If $N(x)$ is a norm on ${\bf R}^n$, then
$N(x - y)$ defines a metric on ${\bf R}^n$.  If $N(x)$ is a
quasinorm on ${\bf R}^n$, then $N(x - y)$ is a quasimetric
on ${\bf R}^n$.  For $\|x\|_p$ with $0 < p < 1$ we have that
$\|x - y\|_p^p$ defines a metric on ${\bf R}^n$.

	For any positive real number $p$ and $x \in {\bf R}^n$
we have that
\begin{equation}
	\|x\|_\infty \le \|x\|_p.
\end{equation}
Conversely we also have that
\begin{equation}
	\|x\|_p \le n^{1/p} \, \|x\|_\infty,
\end{equation}
and in particular $\|x\|_p$ tends to $\|x\|_\infty$ as $p \to \infty$.
More generally, if $p, q$ are positive real numbers with $p < q$, then
\begin{equation}
	\|x\|_q \le \|x\|_p
\end{equation}
for all $x \in {\bf R}^n$, as one can show using the $q = \infty$ case.
One also has that
\begin{equation}
	\|x\|_p \le n^{(1/p) - (1/q)} \, \|x\|_q,
\end{equation}
which can be verified using the convexity of the function $u^{q/p}$
on the nonnegative real numbers.

\section{Lipschitz conditions}
\label{Lipschitz conditions}
\setcounter{equation}{0}

	Let $M$ be a nonempty set, and let $d(x, y)$ be a real-valued
function defined for $x, y \in M$.  We say that $d(x, y)$ is a
\emph{semidistance function} if $d(x, y) \ge 0$ and $d(y, x) = d(x,
y)$ for all $x, y \in M$.  If also $d(x, y)$ satisfies the triangle
inequality, so that
\begin{equation}
	d(x, z) \le d(x, y) + d(y, z)
\end{equation}
for all $x, y, z \in M$, then $d(x, y)$ is said to be a \emph{semimetric}.
A real-valued function $N(x)$ on ${\bf R}^n$ is called a \emph{seminorm}
if $N(x) \ge 0$ for all $x \in {\bf R}^n$, $N(t \, x) = |t| \, N(x)$
for all real numbers $t$ and all $x \in {\bf R}^n$, and
\begin{equation}
	N(x + y) \le N(x) + N(y)
\end{equation}
for all $x, y \in {\bf R}^n$.  If $N(x)$ is a seminorm on ${\bf R}^n$,
then $N(x - y)$ is a semimetric on ${\bf R}^n$.

	Suppose that $M$ is a nonempty set and that $d(x, y)$ is a
semidistance function on $M$.  Define $\rho(x, y)$ for $x, y \in M$
to be the infimum of the sums
\begin{equation}
	\sum_{j=1}^l d(w_j, w_{j-1})
\end{equation}
over all finite sequences $w_1, \ldots, w_l$ of points in $M$ such
that $w_0 = x$ and $w_l = y$.  Thus $\rho(x, y)$ is a semidistance
function on $M$ such that
\begin{equation}
	\rho(x, y) \le d(x, y)
\end{equation}
for all $x, y \in M$, and in fact $\rho(x, y)$ is a semimetric because
it satisfies the triangle inequality by construction.  Conversely if
$\tau(x, y)$ is a semimetric on $M$ such that $\tau(x, y) \le d(x, y)$
for all $x, y \in M$, then $\tau(x, y) \le \rho(x, y)$ for all $x, y
\in M$.

	Let $M_1$, $M_2$ be nonempty sets with semidistance functions
$d_1$, $d_2$, respectively.  If $f$ is a mapping from $M_1$ to $M_2$
and $L$ is a nonnegative real number, then we say that $f$ is
\emph{$L$-Lipschitz} with respect to these semidistance functions if
\begin{equation}
	d_2(f(x), f(y)) \le L \, d_1(x, y)
\end{equation}
for all $x, y \in M_1$.

	Notice that for any mapping $f : M_1 \to M_2$ we have that
\begin{equation}
\label{d_2(f(x), f(y))}
	d_2(f(x), f(y))
\end{equation}
defines a semidistance function on $M_1$, since $d_2$ is a
semidistance function on $M_2$.  If $d_2$ is actually a semimetric on
$M_2$, then (\ref{d_2(f(x), f(y))}) is a semimetric on $M_1$.  Of
course the condition that $f$ be $L$-Lipschitz is the same as saying
that the semidistance function (\ref{d_2(f(x), f(y))}) is less than or
equal to $L$ times $d_1(x, y)$.

	Suppose that $M_3$ is another nonempty set equipped with a
semidistance function $d_3$, and that $f_1$, $f_2$ are mappings from
$M_1$ to $M_2$ and from $M_2$ to $M_3$ which are $L_1$,
$L_2$-Lipschitz for some $L_1, L_2 \ge 0$, respectively.
Then the composition $f_2 \circ f_1$, which is the mapping from
$M_1$ to $M_3$ defined by
\begin{equation}
	(f_2 \circ f_1)(x) = f_2(f_1(x)),
\end{equation}
is $L_1 \, L_2$ Lipschitz, i.e., one multiplies the Lipschitz
constants.

	As a special case, let $M$ be a nonemtpy set with a
semidistance function $d(x, y)$, and let $f$ be a real-valued function
on $M$.  We say that $f$ is $L$-Lipschitz for some $L \ge 0$ if
\begin{equation}
	|f(x) - f(y)| \le L \, d(x, y)
\end{equation}
for all $x, y \in M$.  This is the same as saying that $f$ is
$L$-Lipschitz as a mapping from $M$ into the real line equipped with
the standard metric $|s - t|$, $s, t \in {\bf R}$.

	For a real-valued function $f$ on $M$ the property of being
$L$-Lipschitz is equivalent to
\begin{equation}
	f(x) \le f(y) + L \, d(x, y)
\end{equation}
for all $x, y \in M$.  In particular, if $d$ is a semimetric and
$p$ is a point in $M$, then the function $f_p(x) = d(x, p)$ is
automatically $1$-Lipschitz.

	By contrast, suppose that we take $M$ to be the real line
equipped with the quasimetric $|x - y|^q$ where $q > 1$, and that $f$
is a real-valued Lipschitz function with respect to this quasimetric
on the domain, which means that
\begin{equation}
	|f(x) - f(y)| \le L \, |x - y|^q
\end{equation}
for some $L \ge 0$ and all $x, y \in {\bf R}$.  In this event one
can check that $f$ must be a constant function.

	Let $M$ be a nonempty set with semidistance function $d$, and
suppose that ${\bf R}^n$ is equipped with a seminorm $N$.  Let $f$ be
an $L$-Lipschitz mapping from $M$ to ${\bf R}^n$ with respect to the
semimetric associated to this seminorm, which is to say that
\begin{equation}
	N(f(x) - f(y)) \le L \, d(x, y)
\end{equation}
for all $x, y \in M$.  If $t$ is a real number, then $t \, f(x)$ is
$|t| \, L$-Lipschitz as a mapping from $M$ to ${\bf R}^n$.  If $f_1$,
$f_2$ are $L_1$, $L_2$-Lipschitz mappings from $M$ to ${\bf R}^n$,
then the sum $f_1 + f_2$ is $(L_1 + L_2)$-Lipschitz.

	Now let us specialize to the case of real-valued functions,
and the usual absolute values on the real line.  Specifically, suppose
that $f_1$, $f_2$ are real-valued functions on $M$ which are $L_1$,
$L_2$-Lipschitz, respectively.  One can check that the maximum and
minimum of $f_1$, $f_2$ are $\max(L_1, L_2)$-Lipschitz functions on
$M$.  Under the assumption that $f_1$, $f_2$ are also bounded, the
product of $f_1$ and $f_2$ is Lipschitz too.  More precisely,
if $A_1$, $A_2$ are nonnegative real numbers such that
\begin{equation}
	|f_1(x)| \le A_1, \quad |f_2(x)| \le A_2
\end{equation}
for all $x \in M$, then the product of $f_1$ and $f_2$ is
$(A_1 \, L_2 + A_2 \, L_1)$-Lipschitz.

\section{Connectedness}
\label{connectedness}
\setcounter{equation}{0}

	Let $M$ be a nonempty set equipped with a distance function
$d(x, y)$, and let $\epsilon$ be a positive real number.  By an
\emph{$\epsilon$-chain} in $M$ we mean a finite sequence $w_0, \ldots,
w_l$ of points in $M$ such that
\begin{equation}
	d(w_i, w_{i-1}) \le \epsilon
\end{equation}
for $i = 1, \ldots, l$.  A subset $E$ of $M$ is said to be
\emph{$\epsilon$-connected} if for every pair of points $x, y \in E$
there is an $\epsilon$-chain $w_0, \ldots, w_l$ of points in $E$
such that $w_0 = x$ and $w_l = y$.

	Thus a subset $E$ of $M$ is \emph{not} $\epsilon$-connected if
there is a pair of points in $E$ which cannot be connected by an
$\epsilon$-chain.  This is equivalent to saying that $E$ can be
expressed as the union of two nonempty sets $A$, $B$ such that
\begin{equation}
	d(u, v) > \epsilon
\end{equation}
for all $u \in A$ and $v \in B$.  For in this case there is no
$\epsilon$-chain in $E$ connecting points in $A$ to points in $B$.
Conversely, if $E$ is not chain connected, so that there exist $x, y
\in E$ which cannot be connected by an $\epsilon$-chain of points in
$E$, then one can define $A$ to be the set of points in $E$ which can
be connected by an $\epsilon$-chain of points in $E$ to $x$, and $B$
to be the set of remaining points in $E$, and $A$, $B$ have the
properties described earlier.

	Suppose that $M_1$, $M_2$ are nonempty sets and that $d_1$,
$d_2$ are distance functions on them.  Let $\epsilon > 0$ be given,
and assume that $E$ is an $\epsilon$-connected subset of $M_1$.  Also
let $f$ be an $L$-Lipschitz mapping from $M_1$ to $M_2$.  It is easy
to see that $f(E)$ is then an $L \, \epsilon$-connected subset of
$M_2$.  There are analogous statements for more general ``uniformly
continuous'' mappings, with more complicated adjustments to
$\epsilon$.

	Let $M$ be a nonempty set with distance function $d$ and let
$E$ be a subset of $M$.  We say that $E$ is \emph{chain-connected}
if $E$ is $\epsilon$-connected for all $\epsilon > 0$.  As in the
preceding paragraph, the image of a chain-connected set under a
uniformly continuous mapping is also chain connected.

	A subset $E$ of $M$ is not chain connected if it is not
$\epsilon$-connected for some $\epsilon > 0$, which is equivalent to
saying that there is an $\epsilon > 0$ such that $E$ can be expressed
as the union of two nonempty sets $A$, $B$, with $d(u, v) > \epsilon$
for all $u \in A$ and $v \in B$.  In particular, a subset of a metric
space which is not chain connected is not connected in the sense
discussed in \cite{Rudin}.

	This is equivalent to saying that a subset of a metric space
which is connected in the sense of \cite{Rudin} is also chain
connected.  The converse is not true in general, however.  For
instance, the set of rational numbers, as a subset of the real line
with the standard metric $|x - y|$, is chain connected but not
connected.

	Note that a basic result about connected subsets of metric
spaces is that the image of a connected set under a continuous mapping
is connected too.  For a compact subset of a metric space, one can
show that connectedness and chain connectedness are equivalent.

\section{Hausdorff content}
\label{Hausdorff content}
\setcounter{equation}{0}

	Let $M$ be a nonempty set and $d(x, y)$ a metric on $M$, so
that $M$ becomes a metric space.  A subset $E$ of $M$ is said to be
\emph{bounded} if the set of real numbers of the form $d(x, y)$ for
$x, y \in E$ has an upper bound.  The least upper bound or supremum of
this set of numbers is called the diameter of $E$ and is denoted
$\diam E$, at least if $E$ is not the empty set.  We can also define
the diameter of the empty set to be equal to $0$.  If $E$ is a bounded
subset of $M$, then the closure $\overline{E}$ of $E$, as discussed in
\cite{Rudin}, is also a bounded subset of $M$, and with the same
diameter as $E$.

	Let $\alpha$ be a positive real number.  If $A_1, \ldots, A_k$
is a finite collection of bounded subsets of $M$, then consider the
quantity
\begin{equation}
\label{sum_{i=1}^k (diam A_i)^alpha}
	\sum_{i=1}^k (\diam A_i)^\alpha,
\end{equation}
which is a kind of $\alpha$-dimensional measurement of size of the
$A_i$'s.  If $E$ is a bounded subset of $M$, then we define
$\mathcal{H}^\alpha(E)$, a version of the $\alpha$-dimensional
Hausdorff content of $E$, to be the infimum of the sums
(\ref{sum_{i=1}^k (diam A_i)^alpha}) over all finite collections
$A_1, \ldots, A_k$ of bounded subsets of $M$ such that
\begin{equation}
	E \subseteq \bigcup_{i=1}^k A_i.
\end{equation}

	For instance, we could take $k = 1$ and $A_1 = E$, in which
case (\ref{sum_{i=1}^k (diam A_i)^alpha}) reduces to $(\diam
E)^\alpha$.  In particular we automatically have that
\begin{equation}
	\mathcal{H}^\alpha(E) \le (\diam E)^\alpha.
\end{equation}
As a special case we have that $\mathcal{H}^\alpha$ of the empty set
is equal to $0$.  Also, if $E_1$, $E_2$ are two bounded subsets of $M$
such that $E_1 \subseteq E_2$, then
\begin{equation}
	\mathcal{H}^\alpha(E_1) \le \mathcal{H}^\alpha(E_2),
\end{equation}
because any covering of $E_2$ by finitely many bounded subsets
$A_1, \ldots, A_k$ of $M$ is also a covering of $E_1$.

	For any bounded subset $E$ of $M$,
\begin{equation}
	\mathcal{H}^\alpha(\overline{E}) = \mathcal{H}^\alpha(E).
\end{equation}
Indeed, $\mathcal{H}^\alpha(E) \le \mathcal{H}^\alpha(\overline{E})$
simply because $E \subseteq \overline{E}$.  Conversely, suppose that
$A_1, \ldots, A_k$ are finitely many bounded subsets of $M$ such that
$E$ is contained in the union of the $A_i$'s.  Then the closure of $E$
is contained in the union of the closure of the $A_i$'s, and since the
diameter of the closure of a set is equal to the diameter of a set, it
follows that $\mathcal{H}^\alpha(\overline{E}) \le
\mathcal{H}^\alpha(E)$.

	If $E$, $F$ are bounded subsets of $M$, then
\begin{equation}
\label{subadditivity}
	\mathcal{H}^\alpha(E \cup F) 
		\le \mathcal{H}^\alpha(E) + \mathcal{H}^\alpha(F).
\end{equation}
To see this, let $A_1, \ldots, A_k$ and $B_1, \ldots, B_l$ be arbitrary
finite collections of bounded subsets of $M$ such that
\begin{equation}
	E \subseteq \bigcup_{i=1}^k A_i,
		\quad F \subseteq \bigcup_{j=1}^l B_j.
\end{equation}
Then the combined collection $A_1, \ldots, A_k, B_1, \ldots, B_l$ is a
finite collection of bounded subsets of $M$ such that the union of $E$
and $F$ is contained in the union of the $A_i$'s and $B_j$'s together,
and it follows that
\begin{equation}
	\mathcal{H}^\alpha(E \cup F)
		\le \sum_{i=1}^k (\diam A_i)^\alpha
			+ \sum_{j=1}^l (\diam B_j)^\alpha.
\end{equation}
Because $A_1, \ldots, A_k$ and $B_1, \ldots, B_l$ are arbitrary
finite collections of bounded subsets of $M$ which cover $E$ and $F$,
respectively, we get (\ref{subadditivity}), as desired.

	Now suppose that $M_1$, $M_2$ are metric spaces with metrics
$d_1$, $d_2$, respectively, and that $f$ is an $L$-Lipschitz mapping
from $M_1$ to $M_2$.  If $E$ is a bounded subset of $M_1$, then
\begin{equation}
\label{mathcal{H}^alpha_{M_2}(f(E)) le L^alpha mathcal{H}^alpha_{M_1}(E)}
	\mathcal{H}^\alpha_{M_2}(f(E)) 
		\le L^\alpha \, \mathcal{H}^\alpha_{M_1}(E),
\end{equation}
where the subscripts $M_1$, $M_2$ for $\mathcal{H}^\alpha$ are used to
indicate explicitly in which metric space one is working.  Indeed, if
$A$ is any bounded subset of $M_1$, then $f(A)$ is a bounded subset of
$M_2$, and the diameter of $f(A)$ in $M_2$ is less than or equal to
$L$ times the diameter of $A$ in $M_1$.  If $A_1, \ldots, A_k$ is a
finite collection of bounded subsets of $M_1$ such that $E$ is
contained in the union of the $A_i$'s, then $f(A_1), \ldots, f(A_k)$
is a finite collection of bounded subsets in $M_2$ and $f(E)$ is
contained in the union of the $f(A_i)$'s.  The sum of the $\alpha$
powers of the diameters of the $f(A_i)$'s in $M_2$ is less than or
equal to $L^\alpha$ times the sum of the $\alpha$ powers of the
diameters of the $A_i$'s in $M_1$, which leads to
(\ref{mathcal{H}^alpha_{M_2}(f(E)) le L^alpha
mathcal{H}^alpha_{M_1}(E)}).

\end{document}